\documentclass[12pt,letterpaper]{article}

\usepackage{fullpage}
\usepackage{amssymb}
\usepackage{amsmath}
\usepackage{float}
\usepackage[dvips]{epsfig}
\usepackage{epsf}
\usepackage{alltt}

\def \prend{\vrule depth-1pt height7pt width6pt}
\def \proof{\bigbreak\noindent{\bf Proof.\ \ }}
\def \endpf{{\ \ \prend \medbreak}}

\title{On a Construction of Friedman}

\author{
Jeffrey Shallit\thanks{Research supported in part
by a grant from NSERC.} \ and Ming-wei Wang\\
Department of Computer Science\\
University of Waterloo\\
Waterloo, Ontario, Canada  N2L 3G1\\
{\tt shallit@graceland.uwaterloo.ca} \\
{\tt m2wang@math.uwaterloo.ca} \\
}

\newtheorem{theorem}{Theorem}

\begin{document}

\maketitle

\begin{abstract}
H. Friedman obtained remarkable results about
the longest finite sequence $x$ such that for
all $i \not= j$ the word $x[i..2i]$ is not a subsequence
of $x[j..2j]$.  In this note we consider what happens when
``subsequence'' is replaced by ``subword''.
\end{abstract}

\section{Introduction}
\label{intro}

We say a word $y$ is a {\sl subsequence} of a word $z$ if $y$ can be
obtained by striking out $0$ or more symbols from $z$.  For
example, ``iron'' is a subsequence of ``introduction''.  We say
a word $y$ is a {\sl subword} of a word $z$ if there exist words
$w, x$ such that $z = wyx$.  For example, ``duct'' is a subword
of ``introduction''.\footnote{Europeans sometimes use the
term ``factor'' for what we have called ``subword'', and they use
the term ``subword'' for what we have called ``subsequence''.}

We use the notation $x[k]$ to denote the $k$'th letter chosen
from the string $x$.  (The first letter of a string is $x[1]$.)
We write $x[a..b]$ to denote the subword of $x$ of length $b-a+1$
starting at position $a$ and ending at position $b$.

Recently H. Friedman has found a remarkable construction that
generates extremely large numbers \cite{Friedman:2000a,Friedman:2000b}.
Namely, consider words
over a finite alphabet $\Sigma$ of cardinality $k$.
If an infinite word $\bf x$ has the property that
for all $i, j$ with $0 < i < j$ the subword 
${\bf x}[i..2i]$ is not
a subsequence of ${\bf x}[j..2j]$, call it {\sl self-avoiding}.
We apply the same definition for a finite word $x$ of length
$n$, imposing the additional restriction that $j \leq n/2$.

Friedman shows there are no infinite self-avoiding words over
a finite alphabet.
Furthermore, he shows that
for each $k$ there exists a longest finite self-avoiding word $x$ 
over an alphabet of size $k$.  Call $n(k)$ the length of such a word.
Then clearly $n(1) = 3$ and a simple argument shows
that $n(2) = 11$.  Friedman shows that $n(3)$ is greater
than the incomprehensibly large number $A_{7198} (158386)$,
where $A$ is the Ackermann function.

Jean-Paul Allouche asked what happens when ``subsequence'' is
replaced by ``subword''.   A priori we do not expect results
as strange as Friedman's, since there are no infinite anti-chains
for the partial order defined by ``$x$ is a subsequence of $y$'',
while there {\it are\/} infinite anti-chains for the partial order
defined by ``$x$ is a subword of $y$''.

\section{Main Results}

If an infinite word $\bf x$ has the property that
for all $i, j$ with $0 \leq i < j$ the subword 
${\bf x}[i..2i]$ is not
a subword of ${\bf x}[j..2j]$, we call it {\sl weakly self-avoiding}.
If $x$ is a finite word of length $n$, we apply the same definition
with the additional restriction that $j  \leq n/2$.

\begin{theorem}
Let $ \Sigma  = \lbrace 0, 1, \ldots, k-1 \rbrace$.
\begin{itemize}
	\item[(a)] If $k = 1$,
	the longest weakly self-avoiding word is of length $3$,
	namely $000$.

	\item[(b)] If $k = 2$, there are no weakly self-avoiding
	words of length $> 13$.  There are $8$
	longest weakly self-avoiding words, namely
	$0010111111010$, $0010111111011$,
	$0011110101010$, $0011110101011$
	and the four words obtained by changing $0$ to $1$ and $1$
	to $0$.

	\item[(c)] If $k = 3$, there exists an infinite weakly
	self-avoiding word.
\end{itemize}
\label{main}
\end{theorem}

\proof

(a) If a word $x$ over $\Sigma = \lbrace 0 \rbrace$
is of length $\geq 4$, then it must contain $0000$ as a prefix.
Then $x[1..2] = 00$ is a subword of
$x[2..4] = 000$.

\bigskip

(b)  To prove this result, we create a tree whose root is labeled
with $\epsilon$, the empty word.  If a node's label $x$ is weakly
self-avoiding, then it has two children labeled $x0$ and $x1$.  
This tree is finite if and only if there is a longest 
weakly self-avoiding word.  In this case, the leaves of the tree represent
non-weakly-self-avoiding words that are minimal in the sense that any proper
prefix is weakly self-avoiding.

Now we use a classical breadth-first
tree traversal technique, as follows:
We maintain a queue,
$Q$, and initialize it with the empty word $\epsilon$.  If the queue
is empty, we are done.  Otherwise, we pop the first element $q$ from
the queue and check to see if it is weakly self-avoiding.  If not,
the node is a leaf, and we print it out.  If $q$ is weakly self-avoiding
then we append $q0$ and $q1$ to the end of the queue.

If this algorithm terminates, we have proved that
there is a longest weakly self-avoiding word.  The proof may be
concisely represented by listing the leaves in breadth-first order.
We may shorten the tree by assuming, without loss of generality,
that the root is labeled $0$.

When we perform this procedure, we obtain a tree with 92 leaves,
whose longest label is of length $14$.  The following list 
describes this tree:

\newpage
\begin{alltt}
           0000       00111100       0011010101       001011111011
           0001       00111110       0011010110       001011111100
           0101       00111111       0011010111       001011111110
         001000       01000000       0011101000       001011111111
         001001       01000001       0011101001       001110101000
         001010       01000010       0011101011       001110101001
         001100       01000011       0011110100       001110101010
         010001       01100001       0011110110       001110101011
         010010       01100010       0011110111       001111010100
         010011       01100011       0110000000       001111010110
         011001       01110001       0110000001       001111010111
         011010       01110010       0110000010       011100000000
         011011       01110011       0110000011       011100000001
         011101     0010110100       0111000001       011100000010
         011110     0010110101       0111000010       011100000011
         011111     0010110110       0111000011     00101111110100
       00101100     0010110111     001011110100     00101111110101
       00110100     0010111000     001011110101     00101111110110
       00110110     0010111001     001011110110     00101111110111
       00110111     0010111010     001011110111     00111101010100
       00111000     0010111011     001011111000     00111101010101
       00111001     0010111100     001011111001     00111101010110
       00111011     0011010100     001011111010     00111101010111
\end{alltt}
\centerline{Figure 1:  Leaves of the tree giving a proof of
Theorem~\ref{main} (b)}

\bigskip\bigskip

(c) Consider the word
\begin{eqnarray*}
{\bf x} &=&
2 2 0 1 0 1 1 0 1 1 1 0 1 1 1 1 1 0 1 1 1 1 1 1  1 0 1 1 1 1 1 1 1 1 1 1 1 0
\cdots \\
&=& 2 \, 2 \, 0 \, 1 \, 0 \, 1^2 \, 0 \, 1^3 \, 0 \, 1^5 \, 0 \, 1^7 \,
0 \, 1^{11} \, 0 \, 1^{15} \, 0 \, 1^{23} \, 0 \, 1^{31} \, 0 \, 1^{47} \, 0 \, \cdots
\end{eqnarray*}
where there are $0$'s in positions $3, 5, 8, 12, 18, 26, 38, 54, 78,
110, 158, \ldots$.  More precisely, define
$f_{2n+1} = 5 \cdot 2^n - 2$ for $n \geq 0$, and
$f_{2n} = 7 \cdot 2^{n-1} - 2$ for $n \geq 1$.  Then $\bf x$ has $0$'s
only in the positions given by $f_i$ for $i \geq 1$.

First we claim that if $i \geq 3$, then
any subword of the form ${\bf x}[i..2i]$ contains exactly two $0$'s.
This is easily verified for $i = 3$.  
If $5 \cdot 2^n - 1 \leq i < 7 \cdot 2^n - 1$ and $n \geq 0$,
then there are $0$'s at positions $7 \cdot 2^n - 2$ and
$5 \cdot 2^{n+1} - 2$.  (The next $0$ is at position
$7 \cdot 2^{n+1} - 2$, which is $> 2 (7 \cdot 2^n - 2)$.)
On the other hand, if
$7 \cdot 2^{n-1} - 1 \leq i < 5 \cdot 2^n -1$ for $n \geq 1$, then there
are $0$'s at positions $5 \cdot 2^n - 2$ and
$7 \cdot 2^n -2 $.  (The next $0$ is at position $5 \cdot 2^{n+1} - 2$,
which is $> 2 \cdot (5 \cdot 2^n - 2)$.)

Now we prove that $\bf x$ is weakly self-avoiding.  Clearly
${\bf x}[1..2] = 22$ is not a subword of any subword
of the form ${\bf x}[j..2j]$ for any $j \geq 2$.  Similarly,
${\bf x}[2..4] = 201$ is not a subword of any subword of the
form ${\bf x}[j..2j]$ for any $j \geq 3$.  Now consider
subwords of the form $t := {\bf x}[i..2i]$ and
$t' := {\bf x}[j..2j]$ for $i,j \geq 3$ and $i < j$.  From
above we know $t = 1^u 0 1^v 0 1^w$, and
$t' = 1^{u'} 0 1^{v'} 0 1^{w'}$.  For $t$ to be a subword of $t'$ we must
have $u \leq u'$, $v = v'$, and $w \leq w'$.  But since the blocks of
$1$'s in $\bf x$ are distinct in size, this means that
the middle block of $1$'s in $t$ and $t'$ must occur in the same positions of $\bf x$.
Then $u \leq u'$ implies $i \geq j$, a contradiction.
\endpf

\section{Another construction}

    Friedman also has considered variations on his construction,
such as the following:  let $M_2(n)$ denote the length of the longest
finite word $\bf x$ over $\lbrace 0, 1 \rbrace$
such that ${\bf x}[i..2i]$ is not a subsequence
of ${\bf x}[j..2j]$ for $n \leq i < j$.    We can again consider
this where ``subsequence'' is replaced by ``subword''.

\begin{theorem}
     There exists an infinite word $\bf x$ over $\lbrace 0,1 \rbrace$
such that ${\bf x}[i..2i]$ is not a subword of
${\bf x}[j..2j]$ for all $i, j$ with $2 \leq i < j$.
\end{theorem}

\proof
    Let 
\begin{eqnarray*}
{\bf x} &=& 0 \, 0 \, 1 \, 0 \, 0 \, 1^3 \, 0 \, 1^2 \, 0
\, 1^7 \, 0 \, 1^5 \, 0 \, 1^{15} \, 0 \, 1^{11} \, 0 \, 1^{31}
\, 0 \, 1^{23} \, \cdots \\
&=& 0 \, 0 \, 1 \, 0 \, 0 \, 1^{g_1} \, 0 \, 1^{g_2} \, 0 \, 1^{g_3} \, 0 \,  \cdots  
\end{eqnarray*}
where $g_1 = 3$, $g_2 = 2$, and $g_n = 2g_{n-2} + 1$ for $n \geq 3$.
Then a proof similar to that above shows that every subword
of the form ${\bf x}[i..2i]$ contains exactly two $0$'s, and hence,
since the $g_i$ are all distinct, we have
${\bf x}[i..2i]$ is not a subword of
${\bf x}[j..2j]$ for $j > i > 1$.
\endpf

\end{document}